\documentclass[a4, 10pt]{amsart}
\usepackage{amssymb}
\usepackage{amstext}
\usepackage{amsmath}
\usepackage{amscd}
\usepackage{latexsym}
\usepackage{amsfonts}

\theoremstyle{plain}
\newtheorem{thm}{Theorem}[section]
\newtheorem*{thm*}{Theorem}
\newtheorem*{cor*}{Corollary}

\newtheorem{lem}[thm]{Lemma}

\newtheorem*{claim*}{Claim}

\theoremstyle{definition}
\newtheorem{defn}[thm]{Definition}

\newtheorem{rem}[thm]{Remark}

\theoremstyle{remark}
\newtheorem*{pf}{{\sl Proof}}
\newtheorem*{apf}{{\sl Proof of Theorem \ref{main}(1)}}
\newtheorem*{bpf}{{\sl Proof of Theorem \ref{main}(2)}}

\newtheorem*{ac}{{\sc Acknowledgments}}

\numberwithin{equation}{thm}
\def\Hom{\mathrm{Hom}}
\def\lhom{\underline{\Hom}}
\def\shom{\mathcal{H}{\it om}}

\def\Ext{\mathrm{Ext}}

\def\Im{\mathrm{Im}}

\def\Z{\Bbb Z}

\def\O{\mathcal O}
\def\L{\mathcal L}

\def\div{\mathrm{div}}

\def\Pic{\mathrm{Pic}}

\def\F{{\mathcal F}}
\def\K{{\mathcal K}}

\begin{document}

\title[Totally reflexive modules constructed from curves]{Totally reflexive modules constructed from smooth projective curves of genus $g\geq 2$}
\author{Ryo Takahashi}
\address{Department of Mathematics, School of Science and Technology, Meiji University, 1-1-1 Higashimita, Tama-ku, Kawasaki 214-8571, Japan}
\email{takahasi@math.meiji.ac.jp}
\author{Kei-ichi Watanabe}
\address{Department of Mathematics, College of Humanities and Sciences, Nihon University, 3-25-40 Sakurajosui, Setagaya-Ku, Tokyo 156-8550, Japan}
\email{watanabe@math.chs.nihon-u.ac.jp}
\begin{abstract}
In this paper, from an arbitrary smooth projective curve of genus at least two, we construct a non-Gorenstein Cohen-Macaulay normal domain and a nonfree totally reflexive module over it.
\end{abstract}
\maketitle
\section{Introduction}

About fourty years ago, Auslander \cite{Auslander} introduced a homological invariant for finitely generated modules over a noetherian ring which is called Gorenstein dimension, or G-dimension for short.
He developed the theory of G-dimension with Bridger \cite{AB}.
So far, G-dimension has been studied deeply from various points of view; the details can be found in Christensen's book \cite{Christensen}.

Modules of G-dimension zero are called totally reflexive.
Any finitely generated free module is totally reflexive.
Over a Gorenstein local ring, the totally reflexive modules are precisely the maximal Cohen-Macaulay modules, and hence there exists a nonfree totally reflexive module unless the ring is regular.
Over a non-Gorenstein local ring, on the other hand, it is difficult in general to confirm whether there exists a nonfree totally reflexive module or not.
Avramov and Martsinkovsky \cite[Examples 3.5]{AM} proved that over a Golod local ring that is not a hypersurface (e.g. a Cohen-Macaulay non-Gorenstein local ring with minimal multiplicity \cite[Example 5.2.8]{Avramov2}) every totally reflexive module is free.
Yoshino \cite[Theorem 3.1]{Yoshino} gave several necessary conditions for artinian local rings of Loewy length three to possess nonfree totally reflexive modules.
However, for other cases, we do not have much information about local rings over which there exist nonfree totally reflexive modules.

The purpose of this paper is to construct a non-Gorenstein Cohen-Macaulay normal domain having a nonfree totally reflexive module from any given smooth projective curve of genus at least two.
Nonfree totally reflexive modules have been constructed in several papers \cite{AGP1, AGP2, AM, GP, Yoshino}.
The examples we will show are reflexive modules whose order is two in the class group.
The construction is quite geometric and 
we think that this is the first occasion that the theory of G-dimension is connected to algebraic geometry.
We refer to \cite{Hartshorne} for the basic facts of the geometry 
of algebraic curves.

The following theorem is the main result of this paper.

\begin{thm}\label{main}
Let $C$ be a smooth projective curve of genus $g$ over an algebraically closed field $k$.
Then the following hold.
\begin{enumerate}
\item[{\rm (1)}]
There exists a Weil divisor $D$ on $C$ of degree $g+1$ such that the associated invertible sheaf $\O _C(D)$ is generated by global sections and satisfies $\dim _k H^0(C, \O _C(D))=2,\ H^1(C, \O _C(D))=0$.
\item[{\rm (2)}]
For a divisor $D$ as in (1), set $R=\bigoplus _{n\geq 0} H^0 (C, \O _C(2nD))$ and $M=\bigoplus _{n\geq 0} H^0 (C, \O _C((2n+1)D))$.
Then $R$ is a standard graded normal $k$-algebra of dimension $2$ and type $g$, and $M$ is a finitely generated graded nonfree totally reflexive $R$-module of rank $1$ which has a minimal graded free resolution
$$
\cdots \overset{A}{\to} R(-n)^{\oplus 2} \overset{A}{\to} R(-n+1)^{\oplus 2} \overset{A}{\to} \cdots \overset{A}{\to} R(-1)^{\oplus 2} \overset{A}{\to} R \to M \to 0,
$$
where $A$ is a $2\times 2$ matrix whose entries are elements of $R$ of degree $1$.
\end{enumerate}
\end{thm}

We should note that in the above theorem $R$ is not Gorenstein whenever $g$ is more than or equal to $2$.

In the next section, we will give the definition of a totally reflexive module, and prove the above theorem.

\section{Proof of the theorem}

We start by recalling the definition of a totally reflexive module.

\begin{defn}
Let $R$ be a commutative noetherian ring, and $M$ a finitely generated $R$-module.
We say that $M$ is {\it totally reflexive} if the following two conditions hold:\\
(1) The natural homomorphism $M \to \Hom _R (\Hom _R(M,R),R)$ is an isomorphism, namely, $M$ is reflexive.\\
(2) $\Ext _R ^i(M,R)=\Ext _R^i(\Hom _R(M,R),R)=0$ for every $i>0$.
\end{defn}

Here, we establish the notation which will be used throughout the rest of this section.
Let $k$ be an algebraically closed field, and $C$ a smooth projective curve of genus $g$ over $k$.
Let $\O _C$ denote the structure sheaf of $C$, $\omega _C$ the canonical sheaf of $C$, and $\K _C$ a canonical divisor on $C$.
For a coherent sheaf $\F$ on $C$ and an integer $i$, set $H^i (\F )=H^i (C, \F )$ and $h^i(\F )=h^i(C, \F)=\dim _k H^i(\F )$.
We denote by $\chi (\F )$ the Euler characteristic of $\F$, i.e., $\chi (\F )=h^0(\F )-h^1(\F )$.
For a closed point $P$ on $C$, $k_P$ is the sheaf defined by
$$
\Gamma (U,k_P)=
\begin{cases}
k & \text{if }P\in U,\\
0 & \text{if }P\notin U
\end{cases}
$$
for each open subset $U$ of $C$.

Now, let us prove the first assertion of our theorem.

\begin{apf}
Let $D$ be a divisor on $C$ of degree $g+1$, and let $\L = \O _C  (D)$
 be the invertible sheaf associated to $D$.
Fix a closed point $P$ on $C$.
We have an isomorphism $H^1(\L (-P))\cong H^0(\shom _{\O _C }(\L (-P), 
\omega _C))$ by Serre duality.
Since $\shom _{\O _C }(\L (-P), \omega _C)=\shom _{\O _C }(\O _C (D-P), 
\O _C (\K _C))\cong\O _C (\K _C-D+P)$, the cohomology $H^1(\L(-P))$ does 
not vanish if and only if $\div (x)+(\K _C-D+P)\geq 0$ for some nonzero
 element $x$ in $k(C)$.
 
Now, suppose $H^1(\L (-P))\ne 0$.
Let $E=\div (x)+(\K _C-D+P)\ge 0$ with $0\ne x\in k(C)$, and write $E=P_1+P_2+\cdots +P_{g-2}$, where each $P_i$ is a closed point.
Then $D$ is linearly equivalent to $\K _C+P-(P_1+P_2+\cdots +P_{g-2})$, where $P,P_1,P_2,\dots,P_{g-2}$ run over all closed points.
Let $\Pic_{g+1}(C)$ be the linear equivalence class of the divisors of degree $g+1$ on $C$, and let $\phi : C^{g-1} \to \Pic_{g+1}(C)$ be the morphism of algebraic varieties sending $(P,P_1,P_2,\dots ,P_{g-2})$ to the linear equivalence class of $\K _C+P-(P_1+P_2+\cdots +P_{g-2})$.
Since $\dim C^{g-1}=g-1$ and $\dim \Pic_{g+1}(C)=g$, the morphism $\phi$ cannot be surjective.
Thus, if we take $D$ from outside of the image of $\phi$, then $H^1(\L (-P))=0$ for all closed points $P$, where $\L = \O _C (D)$.

For a closed point $P$, the mapping $H^1(\L (-P)) \to H^1(\L )$ induced from the inclusion $\L (-P) \to \L$ is surjective and we have $H^1(\L)=0$, too.
The Riemann-Roch theorem shows that $h^0(\L)=\chi (\L)=\chi (\O _C)+\deg D=2$.
On the other hand, from the exact sequence $0 \to \L (-P) \to \L \to k_P\otimes _{\O _C}\L \to 0$, we get an exact sequence $H^0(\L) \to H^0(k_P\otimes _{\O _C} \L) \to H^1(\L(-P))=0$.
This induces a surjective homomorphism $H^0(\L)\otimes _k \O _{C,P}\to H^0(k_P\otimes _{\O _C} \L)\otimes _k \O _{C,P}\cong\L _P$, which says that $\L$ is generated by global sections.
Thus, the proof of the first assertion of our theorem is completed.
\qed
\end{apf}

Next, we consider the second assertion of our theorem.
Put $R=\bigoplus _{n\geq 0} H^0(\L ^{\otimes 2n})$ and $M=\bigoplus _{n\geq 0} H^0(\L ^{\otimes 2n+1})$.
Then $R$ is a $2$-dimensional graded normal domain, and $M$ is a finitely generated graded $R$-module; we have $R_n\cdot M_m\subseteq M_{n+m}$ and $M_n\cdot M_m\subseteq R_{n+m+1}$ for any integers $n,m$.
Since $\L$ is generated by global sections and $h^0(\L )=2$, we have an exact sequence
\begin{equation}\label{inv}
0 \to \L ^{-1} \to \O _C^{\oplus 2} \overset{(\alpha,\beta)}{\longrightarrow} \L \to 0
\end{equation}
for some $\alpha,\beta\in H^0(\L)$.

\begin{lem}\label{cpt}
\begin{enumerate}
\item[{\rm (1)}]
For any divisor $E$ on $C$ of degree at least $2g-1$, one has $H^1 ( \O _C (E))=0$.
\item[{\rm (2)}]
There are exact sequences
$$
\begin{cases}
0 \to M_{n-1} \to R_n^{\oplus 2} \overset{(\alpha,\beta)}{\longrightarrow} M_n \to 0,\\
0 \to R_m \to M_m^{\oplus 2} \overset{(\alpha,\beta)}{\longrightarrow} R_{m+1} \to 0
\end{cases}
$$
for $n\geq 0$ and $m\geq 1$.
\item[{\rm (3)}]
One has $\Hom _R(M,R)\cong M(-1)$.
\end{enumerate}
\end{lem}

\begin{pf}
(1) It follows from the Serre duality theorem that there is an isomorphism $H^1(\O _C(E))\cong H^0(\O _C(\K _C-E))$.
Hence if $H^1(\O _C(E))\neq 0$, then $\deg (\K _C-E)\geq 0$, and $\deg E\leq \deg \K _C=2g-2$.

(2) From \eqref{inv} we obtain an exact sequence $0 \to \L^{\otimes 2n-1} \to (\L^{\otimes 2n})^{\oplus 2} \overset{(\alpha,\beta)}{\longrightarrow} \L^{\otimes 2n+1} \to 0$.
This induces an exact sequence
$$
0 \to H^0(\L^{\otimes 2n-1}) \to H^0(\L^{\otimes 2n})^{\oplus 2} \overset{(\alpha,\beta)}{\longrightarrow} H^0(\L^{\otimes 2n+1}) \to H^1(\L^{\otimes 2n-1}).
$$

We claim that $\eta : H^0(\L^{\otimes 2n})^{\oplus 2} \overset{(\alpha,\beta)}{\longrightarrow} H^0(\L^{\otimes 2n+1})$ is a surjective homomorphism.
Indeed, if $n=0$, then $\eta$ is an isomorphism since $\L$ is generated by global sections.
If $n=1$, then $H^1(\L^{\otimes 2n-1})=H^1(\L)=0$.
If $n\geq 2$, then $\deg ((2n-1)D) = (2n-1)(g+1)\geq 2g-1$, hence $H^1(\L^{\otimes 2n-1})=0$ by (1).
Therefore, in any case the homomorphism $\eta$ is surjective, and we obtain an exact sequence $0 \to M_{n-1} \to (R_n)^{\oplus 2} \overset{(\alpha,\beta)}{\longrightarrow} M_n \to 0$.

On the other hand, from \eqref{inv} we get another exact sequence $0 \to
 \L^{\otimes 2m} \to (\L^{\otimes 2m+1})^{\oplus 2} \overset{(\alpha,\beta)}
{\longrightarrow} \L^{\otimes 2m+2} \to 0$, which induces an exact sequence
$$
0 \to H^0(\L^{\otimes 2m}) \to H^0(\L^{\otimes 2m+1})^{\oplus 2}
 \overset{(\alpha,\beta)}{\longrightarrow} H^0(\L^{\otimes 2m+2}) 
\to H^1(\L^{\otimes 2m}).
$$
Since $m\geq 1$, the divisor $2mD$ has degree $2m(g+1)\geq 2g-1$, and $H^1(\L ^{\otimes 2m})=0$ by (1).
Thus we have $0 \to R_m \to M_m^{\oplus 2} \overset{(\alpha,\beta)}{\longrightarrow} R_{m+1} \to 0$.

(3) There are isomorphisms
$$
\Hom _R(M,R) \cong\bigoplus _{n\in\Z}\lhom _R(M,R(n))\cong\bigoplus _{n\in\Z} H^0(\shom _{\O _C}(\widetilde{M}, \widetilde{R(n)})),
$$
and
$$
\shom _{\O _C}(\widetilde{M}, \widetilde{R(n)}) \cong\shom _{\O _C} (\O _C(D), \O _C(2nD))\cong \O _C((2n-1)D).
$$
Therefore we obtain $\Hom _R (M,R)\cong\bigoplus _{n\in\Z} H^0(\O _C((2n-1)D))=M(-1)$.
\qed
\end{pf}

Now we can prove the second assertion of our theorem.

\begin{bpf}
Take the direct sum of copies of the first exact sequence in Lemma \ref{cpt}(2), and we obtain an exact sequence
\begin{equation}\label{fg}
0 \to M(-1) \to R^{\oplus 2} \overset{(\alpha,\beta)}{\longrightarrow} M \to 0
\end{equation}
of graded $R$-modules.
Therefore there is an exact sequence
\begin{equation}\label{cr}
\cdots \to R(-1)^{\oplus 2} \to R^{\oplus 2} \to R(1)^{\oplus 2} \to \cdots
\end{equation}
and the $R$-dual of this sequence is also exact by Lemma \ref{cpt}(3).
It follows that $M$ is a totally reflexive $R$-module.
Noting that $\alpha, \beta$ form a $k$-basis of $H^0(\L )$ and using the exact sequence \eqref{fg}, we easily see that the $R$-module $M$ is nonfree of rank one.

On the other hand, Lemma \ref{cpt}(2) gives a series of surjective homomorphisms
$$
\begin{CD}
R_n @<{(\alpha,\beta)}<< M_{n-1}^{\oplus 2} @<{(\alpha,\beta)}<< R_{n-1}^{\oplus 4} @<{(\alpha,\beta)}<< \cdots @<{(\alpha,\beta)}<< R_1^{\oplus 4^{n-1}}
\end{CD}
$$
for $n\geq 1$.
Note that $\alpha,\beta\in M_0$ and $M_0\cdot M_0\subseteq R_1$.
We easily see that any element of $R_n$ is described as a linear combination of $n$th powers of elements of $R_1$ over $R_0=k$.
This means that $R$ is a standard graded $k$-algebra.

Now, we prove that the ring $R$ has (Cohen-Macaulay) type $g$.
Let $K$ be a canonical module of the Cohen-Macaulay ring $R$.
Then we have $K=\bigoplus _{n\in\Z} H^0(\omega _C\otimes\L ^{\otimes 2n})$.
Since $\omega _C\otimes\L ^{\otimes 2n}\cong\O _C(\K _C+2nD)$ and $\deg (\K _C+2nD)=(2g-2)+2n(g+1)<0$ if $n<0$, we see that $H^0(\omega _C\otimes\L ^{\otimes 2n})=0$ for $n<0$, and hence $K=\bigoplus _{n\geq 0} H^0(\omega _C\otimes\L ^{\otimes 2n})$.

It is enough to show that the $R$-module $K$ is generated by elements of degree zero.
In fact, note that $K_0=H^0(\omega _C)\cong H^1(\O _C)\cong k^{\oplus g}$ by Serre duality.
Hence if the statement is shown, then $K$ is generated by a $k$-basis of $K_0$ as an $R$-module, and we can conclude that the minimal number of generators of the $R$-module $K$ is equal to $g$, equivalently, $R$ has type $g$.

Tensoring $\omega _C\otimes\L ^{\otimes 2n-1}$ with the exact sequence \eqref{inv}, we get an exact sequence $0 \to \omega _C\otimes\L ^{\otimes 2(n-1)} \to (\omega _C\otimes\L ^{\otimes 2n-1})^{\oplus 2} \overset{(\alpha,\beta)}{\longrightarrow} \omega _C\otimes\L ^{\otimes 2n} \to 0$.
From this we obtain an exact sequence
$$
H^0(\omega _C \otimes\L ^{\otimes 2n-1})^{\oplus 2} \overset{(\alpha,\beta)}{\longrightarrow} K_n \to H^1(\omega _C\otimes \L ^{\otimes 2(n-1)}).
$$
Since $\deg (\K _C+2(n-1)D)\geq 2g-1$ for $n\geq 2$, we have $H^1(\omega _C\otimes\L ^{\otimes 2(n-1)})=0$ for $n\geq 2$ by Lemma \ref{cpt}(1), and get a surjective homomorphism $H^0(\omega _C\otimes\L ^{\otimes 2n-1})^{\oplus 2} \overset{(\alpha,\beta)}{\longrightarrow} K_n$ for $n\geq 2$.
Tensoring $\omega _C\otimes \L ^{\otimes 2n-2}$ with \eqref{inv} and making a similar argument, we obtain a surjective homomorphism $K_{n-1}^{\oplus 2} \overset{(\alpha,\beta)}{\longrightarrow} H^0(\omega _C\otimes\L ^{\otimes 2n-1})$ for $n\geq 2$.
Therefore we have
\begin{equation}\label{nn-1}
K_n = \alpha ^2 K_{n-1}+\alpha\beta K_{n-1}+\beta ^2 K_{n-1}\text{ for }n\geq 2.
\end{equation}

The exact sequence \eqref{inv} also gives an exact sequence
\begin{align*}
0 & \to H^0(\omega _C\otimes\L ^{-1}) \to K_0^{\oplus 2} \overset{(\alpha,\beta)}{\longrightarrow} H^0(\omega _C\otimes\L)\\
& \overset{p}{\to} H^1(\omega _C\otimes\L ^{-1}) \overset{q}{\to} H^1(\omega _C)^{\oplus 2} \overset{(\alpha,\beta)}{\longrightarrow} H^1(\omega _C\otimes\L).
\end{align*}
Using the Serre duality theorem, we get $H^0(\omega_C\otimes \L^{-1})\cong H^1(\L)=0$, $H^1(\omega_C\otimes \L)\cong H^0(\L^{-1})=0$ and $H^1(\omega _C\otimes\L ^{-1})\cong H^0(\L)\cong k^{\oplus 2}\cong H^1(\omega _C)^{\oplus 2}$.
Therefore the homomorphism $q$ in the above exact sequence is an isomorphism, 
and thus $p$ is the zero map.
It follows that
\begin{equation}\label{isom}
K_0^{\oplus 2} \overset{(\alpha,\beta)}{\longrightarrow} H^0(\omega _C\otimes\L )\text{ is an isomorphism}.
\end{equation}

Moreover, from \eqref{inv} we obtain the following two exact sequences:
$$
\begin{cases}
0 \to H^0(\O _C) \to K_0^{\oplus 2} \overset{(\alpha,\beta)}{\longrightarrow} R_1 \to H^1(\O _C) \to 0,\\
0 \to K_0 \to H^0(\omega _C\otimes\L)^{\oplus 2} \overset{(\alpha,\beta)}{\longrightarrow} K_1 \to H^1(\omega _C) \to 0.
\end{cases}
$$
Fix a nonzero element $z\in K_0$.
This element induces an exact sequence $0 \to \O _C \overset{z}{\to} \omega _C \to \F \to 0$ with $\dim\F =0$.
There is a commutative diagram
$$
\begin{CD}
0 @>>> H^0(\O _C) @>>> K_0^{\oplus 2} @>{(\alpha,\beta)}>> R_1 @>>> H^1(\O _C) @>>> 0,\\
@. @V{z}VV @V{z}VV @V{z}V{g}V @V{z}V{f}V \\
0 @>>> K_0 @>>> H^0(\omega _C\otimes\L)^{\oplus 2} @>{(\alpha,\beta)}>h> K_1 @>>> H^1(\omega _C) @>>> 0.
\end{CD}
$$
with exact rows.
Since $H^1(\F)=0$, the homomorphism $f$ is surjective.
Diagram chasing shows that $K_1=\Im\,g + \Im\,h$.
We have $\Im\,g=R_1z$, and $\Im\,h=\alpha ^2K_0+\alpha\beta K_0+\beta ^2K_0$ by \eqref{isom}.
Therefore we obtain
\begin{equation}\label{0}
K_1 = \alpha ^2K_0 + \alpha\beta K_0+ \beta ^2 K_0 +R_1 z.
\end{equation}
Note that $\alpha ^2$, $\alpha\beta$ and $\beta ^2$ are elements of $R_1$.
Putting together \eqref{0} and \eqref{nn-1}, we see that any element of $K_n$ can be described as a linear combination of elements of $K_0$ over $R$.
Thus, the $R$-module $K$ is generated by elements of degree zero, and $R$ has type $g$.
See also \cite[Example (2.12)]{Tomari}.
\qed
\end{bpf}

\begin{rem}\label{061006}
(1) The Hilbert series of the graded ring $R$ is as follows:
$$
H_R(t)=\frac{1+(g+1)t+gt^2}{(1-t)^2}.
$$

Indeed, it is seen by Lemma \ref{cpt}(1) that $h^1(\O _C(2nD))=0$ for any $n\geq 1$.
Applying the Riemann-Roch theorem, we get the Hilbert function of $R$:
\begin{align*}
H(R,n) & =\dim _k R_n =h^0(\O _C(2nD))=\chi (\O _C(2nD)) \\
& =\chi (\O _C)+\deg (2nD) =(1-g)+2n(g+1)
\end{align*}
for $n\geq 1$, and $H(R,0)=h^0(\O _C)=1$.
Thus we obtain the Hilbert series of $R$:
$$
H_R(t)=1+\sum _{n\geq 1}((1-g)+2n(g+1))t^n=\frac{1+(g+1)t+gt^2}{(1-t)^2},
$$
as desired.

(2) Using the fact that $R$ admits a totally reflexive module, one can also show that $R$ has type $g$, as follows.

Since $R$ is a standard graded algebra over an infinite field $k$, we can choose a homogeneous system of parameters $x,y$ of $R$ in $R_1$ (cf. \cite[Theorem 1.5.17(c)]{BH}), and the statement (1) yields the Hilbert series of the residue ring $R/(x,y)$:
$$
H_{R/(x,y)}(t)=1+(g+1)t+gt^2.
$$
This shows that the residue ring $R/(x,y)$ is an artinian ring the cube of whose graded maximal ideal is zero.
On the other hand, using \cite[Lemma (1.3.5)]{Christensen}, we see that $M/(x,y)M$ is a nonfree totally reflexive $R/(x,y)$-module.
Hence, according to \cite[Theorem 3.1]{Yoshino}, the ring $R/(x,y)$ has type $g$, and so does the ring $R$.

(3) It is seen from the exact sequence \eqref{cr} that all the Betti numbers of the $R$-module $M$ are equal to two.
This especially says that $M$ has finite complexity, hence $M$ has lower complete intersection dimension zero; see \cite{Gerko} or \cite{Avramov} for the details.

(4) The ring $R/(x,y)$ obtained in (2) is an artinian local ring of Loewy length $3$, just like the ring  $P$ in \cite[Proposition (3.4)]{GP}, and both $R/(x,y)$ and $P$ yield totally reflexive modules with constant Betti numbers of value $2$.
The authors do not know how the two constructions differ.
However, since the Hilbert series of $P$ is $H_P(t)=1+4t+3t^2$, the ring $R/(x,y)$ is not isomorphic to $P$ unless $R$ has genus $g=3$.
\end{rem}


\begin{ac}
The authors would like to thank an anonymous referee for pointing out several references relevant to this paper and suggesting that they should mention Remark \ref{061006}(4).
\end{ac}


\end{document}